\documentstyle[12pt]{article}
\date{}
\author{Valerii Dryuma\thanks{Work supported in part by Grant 14-0100389 and DAAD, Yena, Germany (2013)}\\[5mm]
{\it Institute of Mathematics and Informatics, AS RM,}\\[3mm] {\it
5 Academiei Street, 2028 Kishinev, Moldova},\\[3mm]{\it e-mail:
valdryum@gmail.com }}
\title{On vector field generated by the Hopf map $S^3$ on $S^2$}
\newtheorem{theorem}{Theorem}

\begin{document}
\maketitle
\date{}
\maketitle
\begin{abstract}
\ \ \ \ The examples of solutions of the system of differential equations generated by the Hopf map $S^3\rightarrow S^2$ are constructed. Their properties are discussed.
\end{abstract}


\medskip

\section{Introduction}

      The subject of our consideration is the system  of ODE`s in $E^3$ -space [1]
      \begin{equation} \label{dr:eq01}
 {\frac {d}{dt}}x(t)=8\,{\frac {4\,zx-y\left ({x}^{2}+{y}^{2}+{z}^{2}
\right )+4\,y}{\left ({x}^{2}+{y}^{2}+{z}^{2}+4\right )^{2}}},
$$$${\frac {d}{dt}}y(t)=8\,{\frac {4\,zy+x\left ({x}^{2}+{y}^{2}+{z}^{2}
\right )-4\,x}{\left ({x}^{2}+{y}^{2}+{z}^{2}+4\right )^{2}}},
$$$${\frac {d}{dt}}z(t)={\frac {24\,{x}^{2}+24\,{y}^{2}-8\,{z}^{2}-\left (
{x}^{2}+{y}^{2}+{z}^{2}\right )^{2}-16}{\left ({x}^{2}+{y}^{2}+{z}^{2}
+4\right )^{2}}}
\end{equation}
 associated with the Hopf map $S^3\rightarrow S^2$ of three dimensional sphere  with equation $z_1\bar z_1+z_2\bar z_2=1$ on two -dimensional  sphere $S^2$ considered  as set of the points $[\frac{z_1}{z_2},1]$ of the complex space $C^2$.

   As it was shown in [1] from the system (\ref{dr:eq01}) can be obtained the  relation
      \begin{equation}\label{dr:eq02}
     \frac {d \arctan(y(t)/x(t))}{dt}=8 \frac{x^2+y^2+z^2 -4}{(x^2+y^2+z^2 +4)^2},~$$$$
     \frac {d( x(t)^2+y(t)^2)}{dt}=64 z \frac{x^2+y^2}{(x^2+y^2+z^2 +4)^2}
\end{equation}
and with the help of which was made qualitative analysis of the behavior of integral curves of the system (\ref{dr:eq01}) .

    Here we consider properties of solutions of the system [1] more details.

\section{Transformation}

\begin{theorem}
      In spherical system of coordinates $r,\phi,\psi$
       $$
       x(t)=r(t)\cos(\phi(t))\sin(\psi(t)),~
y(t)=r(t)\sin(\phi(t))\sin(\psi(t)),~
z(t)=r(t)\cos(\psi(t))
$$ the system of equations (\ref{dr:eq01}) takes the form
 \begin{equation}\label{dr:eq03}
   {\frac {d}{dt}}\phi(t)=-8\,{\frac {\left (r(t)\right )^{2}-4}{16+
\left (r(t)\right )^{4}+8\,\left (r(t)\right )^{2}}},$$$$
{\frac {d}{dt}}r(t)={\frac {\left (\left (r(t)\right )^{4}+8\,\left (r
(t)\right )^{2}-64\,\left (r(t)\right )^{2}\left (\sin(\psi(t))\right
)^{2}+16\right )\cos(\psi(t))}{16+\left (r(t)\right )^{4}+8\,\left (r(
t)\right )^{2}}},$$$$
{\frac {d}{dt}}\psi(t)=-{\frac {\sin(\psi(t))\left (\left (r(t)\right
)^{4}+40\,\left (r(t)\right )^{2}-64\,\left (r(t)\right )^{2}\left (
\sin(\psi(t))\right )^{2}+16\right )}{r(t)\left (16+\left (r(t)\right
)^{4}+8\,\left (r(t)\right )^{2}\right )}}.
\end{equation}
\end{theorem}

       Now from the system (\ref{dr:eq03})  we obtain the equation
   \begin{equation}\label{dr:eq53}
   \left ({\frac {d}{dr}}\psi(r)\right ){r}^{5}\cos(\psi(r))-56\,\left ({
\frac {d}{dr}}\psi(r)\right ){r}^{3}\cos(\psi(r))+64\,\left ({\frac {d
}{dr}}\psi(r)\right ){r}^{3}\left (\cos(\psi(r))\right )^{3}+$$$$+16\,
\left ({\frac {d}{dr}}\psi(r)\right )r\cos(\psi(r))+\sin(\psi(r)){r}^{
4}-24\,\sin(\psi(r)){r}^{2}+$$$$+64\,\sin(\psi(r)){r}^{2}\left (\cos(\psi(r
))\right )^{2}+16\,\sin(\psi(r))=0,
 \end{equation}
which  after substitution
\begin{equation}\label{dr:eq04}
\psi(r)=\arcsin(\sqrt {H(r)})
\end{equation}
 takes the form
\begin{equation}\label{dr:eq05}
\left ({r}^{5}+8\,{r}^{3}-64\,{r}^{3}H(r)+16\,r\right ){\frac {d}{dr}}
H(r)+2\,H(r){r}^{4}+80\,H(r){r}^{2}-128\,\left (H(r)\right )^{2}{r}^{2
}+32\,H(r)=0.
\end{equation}

   General solution of the equation (\ref{dr:eq05}) is determined from implicit relation
\begin{equation}\label{dr:eq06}
4\,{\it C_1}\,{\it BesselI}(0,1/2\,\sqrt {H(r)}r)+{\it C_1}\,{\it
BesselI}(0,1/2\,\sqrt {H(r)}r){r}^{2}-$$$$-8\,{\it C_1}\,{\it BesselI}(1,1
/2\,\sqrt {H(r)}r)\sqrt {H(r)}r+4\,{\it BesselK}(0,-1/2\,\sqrt {H(r)}r
)+$$$$+{\it BesselK}(0,-1/2\,\sqrt {H(r)}r){r}^{2}-8\,{\it BesselK}(1,-1/2
\,\sqrt {H(r)}r)\sqrt {H(r)}r=0,
\end{equation}
where $\it C_1$ is arbitrary.

          In result we obtain that the function $\psi(r)$  from (\ref{dr:eq04}) can be expressed through the function $H(r)$ and then  with a help of the first equation  of the system (\ref{dr:eq03}) can be studied properties of the function $\phi(r)$.

\section{Associated p.d.e.}

   Let us consider   partial differential equation of first order which determine first integral $\Phi(r,\phi,\psi)=C$ of the system (\ref{dr:eq03})
\begin{equation}\label{dr:eq07}
 \left (-\sin(\psi){r}^{4}+24\,\sin(\psi){r}^{2}-64\,\sin(\psi){r}^{2}
\left (\cos(\psi)\right )^{2}-16\,\sin(\psi)\right ){\frac {\partial }
{\partial \psi}}\Phi(r,\phi,\psi)+$$$$\left (-8\,{r}^{3}+32\,r\right ){
\frac {\partial }{\partial \phi}}\Phi(r,\phi,\psi)+$$$$+\left (\cos(\psi){r
}^{5}-56\,\cos(\psi){r}^{3}+64\,\left (\cos(\psi)\right )^{3}{r}^{3}+
16\,\cos(\psi)r\right ){\frac {\partial }{\partial r}}\Phi(r,\phi,\psi
)=0.
\end{equation}

     After substitution of the form
     $$
     \psi=\arcsin(\xi),~
\Phi(r,\phi,\psi)=E(r,\phi,\xi),$$$$
{\frac {\partial }{\partial \psi}}\Phi(r,\phi,\psi)=\left ({\frac {
\partial }{\partial \xi}}E(r,\phi,\xi)\right )\cos(\psi)
$$
and taking in consideration the condition
$$
\cos(\psi)=\sqrt {1-{\xi}^{2}}
$$
we obtain from (\ref{dr:eq07})  the equation
\begin{equation}\label{dr:eq08}
\left ({\frac {\partial }{\partial r}}E(r,\phi,\xi)\right )\sqrt {1-{
\xi}^{2}}{r}^{5}+8\,\left ({\frac {\partial }{\partial r}}E(r,\phi,\xi
)\right )\sqrt {1-{\xi}^{2}}{r}^{3}-64\,\left ({\frac {\partial }{
\partial r}}E(r,\phi,\xi)\right )\sqrt {1-{\xi}^{2}}{r}^{3}{\xi}^{2}+$$$$+
16\,\left ({\frac {\partial }{\partial r}}E(r,\phi,\xi)\right )\sqrt {
1-{\xi}^{2}}r-8\,\left ({\frac {\partial }{\partial \phi}}E(r,\phi,\xi
)\right ){r}^{3}+32\,\left ({\frac {\partial }{\partial \phi}}E(r,\phi
,\xi)\right )r-$$$$-\left ({\frac {\partial }{\partial \xi}}E(r,\phi,\xi)
\right )\sqrt {1-{\xi}^{2}}\xi\,{r}^{4}-40\,\left ({\frac {\partial }{
\partial \xi}}E(r,\phi,\xi)\right )\sqrt {1-{\xi}^{2}}\xi\,{r}^{2}+$$$$+64
\,\left ({\frac {\partial }{\partial \xi}}E(r,\phi,\xi)\right )\sqrt {
1-{\xi}^{2}}{\xi}^{3}{r}^{2}-16\,\left ({\frac {\partial }{\partial
\xi}}E(r,\phi,\xi)\right )\sqrt {1-{\xi}^{2}}\xi=0.
\end{equation}

   It has solution of the form $\Phi(r,\phi,\psi)={\it F_2}(\phi)+H(r,\psi)
$,
where
$$
{\frac {d}{d\phi}}{\it F_2}(\phi)={\it c}_{{2}}
$$
and the function $H(r,\psi)$ satisfies to the equation
\begin{equation}\label{dr:eq09}
\left ({\frac {\partial }{\partial r}}H(r,\psi)\right )\cos(\psi){r}^{
5}-8\,\left ({\frac {\partial }{\partial r}}H(r,\psi)\right )\cos(\psi
){r}^{3}+16\,\left ({\frac {\partial }{\partial r}}H(r,\psi)\right ){r
}^{3}\cos(3\,\psi)+$$$$+16\,\left ({\frac {\partial }{\partial r}}H(r,\psi)
\right )\cos(\psi)r-8\,{\it c}_{{2}}{r}^{3}+32\,{\it c}_{{2}}r-
\left ({\frac {\partial }{\partial \psi}}H(r,\psi)\right )\sin(\psi){r
}^{4}+8\,\left ({\frac {\partial }{\partial \psi}}H(r,\psi)\right )
\sin(\psi){r}^{2}-$$$$-16\,\left ({\frac {\partial }{\partial \psi}}H(r,
\psi)\right ){r}^{2}\sin(3\,\psi)-16\,\left ({\frac {\partial }{
\partial \psi}}H(r,\psi)\right )\sin(\psi)=0
\end{equation}
where ${\it c_{{2}}}$ is parameter.

\subsection{The $(u,v)$-method of integration}

 To derive particular solutions of the partial nonlinear
       differential equation
\begin{equation}\label{Dr10}
F(x,y,f_x,f_y,f_{xx},f_{xy},f_{yy},f_{xxx},f_{xyy},f_{xxy},..)=0
 \end{equation}
    we  use  parametric presentation of the functions and variables
\begin{equation}\label{Dr11}
f(x,y)\rightarrow u(x,t),\quad y \rightarrow v(x,t),\quad
f_x\rightarrow u_x-\frac{u_t}{v_t}v_x,~\quad f_y \rightarrow \frac{u_t}{v_t},$$
$$ f_{yy} \rightarrow \frac{(\frac{u_t}{v_t})_t}{v_t},~
f_{xy} \rightarrow \frac{(u_x-\frac{u_t}{v_t}v_x)_t}{v_t},~ f_{xx} \rightarrow (f_x)_x-\frac{(f_x)_t}{v_t}v_x
\end{equation}
where variable $t$ is considered as parameter [2].

  Remark that conditions of compatibility   $ f_{xy}=f_{yx}$
are fulfilled at the such type of presentation.

  In result instead of equation (\ref{Dr10}) one gets the
  relation between the new variables $u(x,t)$, $v(x,t)$ and
  their partial derivatives
\begin{equation}\label{Dr12}
\Psi(u,v,u_x,u_t,v_x,v_t...)=0.
  \end{equation}

   This value coincides with the initial equation under the condition that $v(x,t)=t$
and takes a more general view in the choice of $u(x,t)=F(\omega(x,t),\omega_t...)$
and $v(x,t)=\Phi(\omega(x,t),\omega_t...)$.

    In some cases the integration of equation (\ref{Dr12}) is a more simple
    problem than the equation (\ref{Dr10}).

     With help of solutions of equation (\ref{Dr12}) $u(x,t),~v(x,t)$ the solutions of initial equation (\ref{Dr10}) $f(x,y)$ are
obtained by exclusion of the parameter $t$ from the conditions
\begin{equation}\label{Dr14}
f(x,y)=u(x,t),~y=v(x,t).
\end{equation}

Next, we apply this method for integrating equations (\ref{dr:eq09}).

    In result of transformation of the function $H(r,\psi)$ and its derivatives
      with accordance of the rules (\ref{Dr11})
           $$
     H( \xi,\psi) =u \left( \xi,\eta \right),~
r=v( \xi,\psi)
 $$
\begin{equation}\label{dr:eq100}
{\frac {\partial }{\partial r}}H(r,\psi)={\frac {{\frac {\partial }{
\partial \xi}}u(\xi,\psi)}{{\frac {\partial }{\partial \xi}}v(\xi,\psi
)}},~~~
{\frac {\partial }{\partial \psi}}H(r,\psi)={\frac {\partial }{
\partial \psi}}u(\xi,\psi)-{\frac {\left ({\frac {\partial }{\partial
\psi}}v(\xi,\psi)\right ){\frac {\partial }{\partial \xi}}u(\xi,\psi)}
{{\frac {\partial }{\partial \xi}}v(\xi,\psi)}}
\end{equation}
the equation (\ref{dr:eq09}) takes form of the relation between the functions $u(\xi,\psi)$ and $v(xi,\psi)$
\begin{equation}\label{dr:eq10}
 \left (8\,\sin(\psi){\xi}^{2}-\sin(\psi){\xi}^{4}-
16\,{\xi}^{2}\sin(3\,\psi)
-16\,\sin(\psi)\right ){
\frac {\partial }{\partial \psi}}u(\xi,\psi){\frac {\partial }{\partial \xi}}v(\xi,
\psi)+$$$$+\left (\sin(\psi){\xi}^{4
}-8\,\sin(\psi){\xi}^{2}+16\,{\xi}^{2}\sin(3\,\psi)+16\,\sin(\psi)
\right )\left ({\frac {\partial }{\partial \xi}}u(\xi,\psi)\right ){
\frac {\partial }{\partial \psi}}v(\xi,\psi)+$$$$+\left (16\,\cos(\psi)\xi-
8\,\cos(\psi){\xi}^{3}+\cos(\psi){\xi}^{5}+16\,{\xi}^{3}\cos(3\,\psi)
\right ){\frac {\partial }{\partial \xi}}u(\xi,\psi)-$$$$-8\,{\it c}_{{2}
}{\xi}^{3}{\frac {\partial }{\partial \xi}}v(\xi,\psi)+32\,{\it c}_{
{2}}\xi\,{\frac {\partial }{\partial \xi}}v(\xi,\psi) =0,
\end{equation}
 where the variable $\xi$ is considered as parameter.

 From here we can obtain the p.d.e by  various ways.
     One of them can be obtained after the substitution of the form
\begin{equation}\label{dr:eq101}
v(\xi,\psi)=\xi\,{\frac {\partial }{\partial \xi}}\rho(\xi,\psi)-\rho(
\xi,\psi),~~
u(\xi,\psi)={\frac {\partial }{\partial \xi}}\rho(\xi,\psi)
\end{equation}
that leads to the equation
\begin{equation}\label{dr:eq11}
\left (\sin(\psi){\xi}^{4}-8\,\sin(\psi){\xi}^{2}+16\,{\xi}^{2}\sin(3
\,\psi)+16\,\sin(\psi)\right ){\frac {\partial }{\partial \psi}}\rho(
\xi,\psi)+$$$$+{\xi}^{6}\cos(\psi)-8\,{\xi}^{4}\cos(\psi)+16\,{\xi}^{4}\cos
(3\,\psi)+16\,{\xi}^{2}\cos(\psi)-8\,{\it c}_{{2}}{\xi}^{3}+32\,{
\it c}_{{2}}\xi=0
\end{equation}
to the function $\rho(\xi,\psi)$.

     Solution of this equation has the form
\begin{equation}\label{dr:eq12}
\rho(\xi,\psi)={\frac {{\xi}^{5}\ln (\chi)}{16+40\,{\xi}^{2}+{\xi}^{4}
}}+$$$$+16\,{\frac {{\xi}^{3}\ln ({\xi}^{4}{\chi}^{4}+2\,{\xi}^{4}{\chi}^{2
}+{\xi}^{4}+40\,{\xi}^{2}{\chi}^{4}-176\,{\xi}^{2}{\chi}^{2}+40\,{\xi}
^{2}+16\,{\chi}^{4}+32\,{\chi}^{2}+16)}{16+40\,{\xi}^{2}+{\xi}^{4}}}+$$$$+8
\,{\frac {{\xi}^{3}\ln (\chi)}{16+40\,{\xi}^{2}+{\xi}^{4}}}-\xi\,\ln (
{\chi}^{2}+1)+16\,{\frac {\xi\,\ln (\chi)}{16+40\,{\xi}^{2}+{\xi}^{4}}
}-$$$$-64\,{\it \_c}_{{2}}{\xi}^{6}\arctan(1/32\,{\frac {2\,{\chi}^{2}
\left (16+40\,{\xi}^{2}+{\xi}^{4}\right )+32-176\,{\xi}^{2}+2\,{\xi}^{
4}}{\sqrt {16\,{\xi}^{2}-24\,{\xi}^{4}+{\xi}^{6}}}})\times$$$$\times\left (16+40\,{\xi
}^{2}+{\xi}^{4}\right )^{-1}{\frac {1}{\sqrt {16\,{\xi}^{2}-24\,{\xi}^
{4}+{\xi}^{6}}}}-$$$$-8\,{\frac {{\it c}_{{2}}{\xi}^{4}\ln (\chi)}{16+40
\,{\xi}^{2}+{\xi}^{4}}}+$$$$+256\,{\it c}_{{2}}{\xi}^{4}\arctan(1/32\,{
\frac {2\,{\chi}^{2}\left (16+40\,{\xi}^{2}+{\xi}^{4}\right )+32-176\,
{\xi}^{2}+2\,{\xi}^{4}}{\sqrt {16\,{\xi}^{2}-24\,{\xi}^{4}+{\xi}^{6}}}
})\times$$$$\times\left (16+40\,{\xi}^{2}+{\xi}^{4}\right )^{-1}{\frac {1}{\sqrt {16\,
{\xi}^{2}-24\,{\xi}^{4}+{\xi}^{6}}}}+$$$$+32\,{\frac {{\it \_c}_{{2}}{\xi}^
{2}\ln (\chi)}{16+40\,{\xi}^{2}+{\xi}^{4}}}+{\it F_1}(\xi),
\end{equation}
where $\chi=\tan(1/2\,\psi)$ and ${\it F_1}(\xi)$-is arbitrary.

     On base of of solutions of the equation (\ref{dr:eq12}) can be find the functions $u(\xi,\psi)$ and $v(\xi,\psi)$ from (\ref{dr:eq101})
 and then get  solutions of the equation  (\ref{dr:eq09}) in parametric form
$$
H(r,\psi)=u(\xi,\psi),~~ r=v(\xi,\psi).
$$

\centerline{\bf References:}

\smallskip
\noindent 1. Aminov Yu.A., {\it Geometriya vektornogo polya}. M.: Nauka. Gl.red.fiz.-mat.lit.,(1990), p.1--208.

\smallskip
\noindent 2. Dryuma V.S., {\it The Riemann and Einsten-Weyl
geometries in theory of differential equations, their applications
and all that}. A.B.Shabat et all.(eds.), New Trends in
Integrability and Partial Solvability, Kluwer Academic Publishers,
Printed in the Netherlands , 2004, p.115--156.
\end{document}